\documentclass[11pt]{amsart}
\usepackage{amsmath,amssymb,amsfonts,url,mathptmx,tikz,geometry}
\usepackage{color}
\usepackage{appendix}
\numberwithin{equation}{section}

\newtheorem{thm}{Theorem}[section]

\usepackage{amsthm,a0size,bm,indentfirst,wasysym}

\usepackage{amsmath}

\allowdisplaybreaks[4]

\begin{document}
	\title[Compactness and existence of solutions]{A complementary result on a singular mean field equation with a sign-changing potential function}

 \author{Lina Wu}\footnote{Lina Wu is partially supported by the National Natural Science Foundation of China (12201030). }

\address{  Lina Wu\\
	School of Mathematics and Statistics \\
	Beijing Jiaotong University \\
	Beijing, 100044, China }
\email{lnwu@bjtu.edu.cn}

\date{\today}
	

\begin{abstract}
	In this note, we study the singular mean field equation defined on a Riemann surface with a sign-changing potential function. We prove if some singular sources happen to be placed on the zero-level curve of the potential function, a priori estimate can still be obtained. As a consequence of this estimate, existence and multiplicity results can still be obtained based on the topology of the manifold.
\end{abstract}
 
	
	
\maketitle

\section{Introduction}    

Let $(M,g)$ be a compact Riemann surface with Gauss curvature $K_g$, one important mean field type equation with singular sources is of the following form
\begin{equation}\label{main-2}
-\Delta_gv+2K_g(x)=2K(x)e^v-4\pi\sum_{j=1}^m\alpha_j\delta_{p_j}
\end{equation}
where $\alpha_1,...,\alpha_m>-1$ are real numbers, $p_1,...,p_m$ are the locations of singular sources, $K$ is a smooth function related to the change of metric. Equation (\ref{main-2}) has wide and profound applications in conformal geometry and physics. In spite of the vast literature concerning various aspects of this famous equation, very few of them talk about the solutions of (\ref{main-2}) when $K$ changes signs. In this respect we address an important work of De Marchis, et, al \cite{jmpa-4}. It is based on this work that we state our main result. 

In \cite{jmpa-4}, the authors proved a number of existence, and multiplicity results based on a compactness result when the function $K$ changes signs. In order to state our results, we denote $\Gamma$ as zero-level curve of $K$: $\Gamma:=\{x\in M,\; K(x)=0\}$, and list the major assumptions in \cite{jmpa-4}:
\begin{align*}
    &\mbox{(H1)\quad is a sign-changing } C^{2,\alpha} \mbox{ function with } \nabla K(x)\neq 0 \mbox{ for any } x\in \Gamma. \\
    &\mbox{(H2)}\quad p_j\not \in \Gamma \mbox{ for all } j\in \{1,...,m\}.
\end{align*}

The major assumption (H2) is postulated in all the main theorems of \cite{jmpa-4}. The main purpose of this short note is to prove that (H2) can be removed from all of them. For example, one main result in \cite{jmpa-4} is Theorem 2.3, with the removal of (H2), the improved version now states the following:
\begin{thm}
    Let $\alpha_1,...,\alpha_l\in (0,1]$, and $\lambda\in (8\pi,16\pi)\setminus \Lambda$. If (H1), (H4) are satisfied then $(*)_{\lambda}$ admits a solution.
\end{thm}

We refer the readers to \cite{jmpa-4} about the notations in the theorem above. The different part of their proof is to establish the a priori estimate, as Theorem 2.1 in \cite{jmpa-4}, without the assumption (H2). 
With such a priori estimate, all the existence and multiplicity results can be obtained by means of variational methods as in \cite{jmpa-4}. 
In the next section, we put the main difference in the proof on a priori estimate.

\section{A Priori Estimate}
In this section, we aim to complete the proof on a priori estimate under the only assumption (H1). To be specific, we need to establish the a priori estimate when a positive singular source is placed on $\Gamma$ by means of the method of moving planes.

For simplicity we assume that $0$ is a singular source on $\Gamma$ with strength $4\pi\alpha$ ($\alpha>0$). The equation around $0$ can be written as, in local coordinates, 
$$\Delta u+|x|^{2\alpha}W(x)e^{u}=0,\quad \mbox{in}\ \ B(0,\tau),$$
where $B(0,\tau)$ stands for the ball centered at the origin with radius $\tau>0$. Here $W$ is the product of $K$ and some positive smooth functions. Since (H1) holds, we know that $\Gamma$ is a $C^{2,\alpha}$ curve. Through Kelvin transform and a rotation, we can assume that $\Gamma$ is tangent to $x_2$ axis at the origin, and $\Gamma$ is contained in $\{x_1<0\}$. Moreover, (H1) implies $\nabla W\neq 0$ on $\Gamma$. We illustrate with the following figure as in \cite{jmpa-4}.

\bigskip


\begin{center}

    \tikzset{every picture/.style={line width=0.75pt}} 

\begin{tikzpicture}[x=0.75pt,y=0.75pt,yscale=-1,xscale=1]

\draw  [fill={rgb, 255:red, 195; green, 195; blue, 195 }  ,fill opacity=0.65 ][line width=0.75]  (141.54,239.32) .. controls (326.97,189.34) and (326.52,141.06) .. (140.17,94.45) ;
\draw  [fill={rgb, 255:red, 155; green, 155; blue, 155 }  ,fill opacity=1 ][line width=0.75]  (141.59,230.13) .. controls (260.06,186.61) and (259.79,145.34) .. (140.79,106.3) ;
\draw [line width=0.75]  (97,166.73) -- (316.22,166.73)(230.16,74.77) -- (230.16,249.31) (309.22,161.73) -- (316.22,166.73) -- (309.22,171.73) (225.16,81.77) -- (230.16,74.77) -- (235.16,81.77)  ;
\draw    (279.23,75.15) -- (280.74,247.8) ;
\draw    (141.09,78.15) -- (142.6,250.81) ;
\draw  [color={rgb, 255:red, 155; green, 155; blue, 155 }  ,draw opacity=1 ][fill={rgb, 255:red, 198; green, 195; blue, 195 }  ,fill opacity=0.65 ] (345.83,191.12) -- (364.01,191.12) -- (364.01,200.72) -- (345.83,200.72) -- cycle ;
\draw  [color={rgb, 255:red, 128; green, 128; blue, 128 }  ,draw opacity=1 ][fill={rgb, 255:red, 155; green, 155; blue, 155 }  ,fill opacity=1 ] (345.83,217.32) -- (364.01,217.32) -- (364.01,226.92) -- (345.83,226.92) -- cycle ;

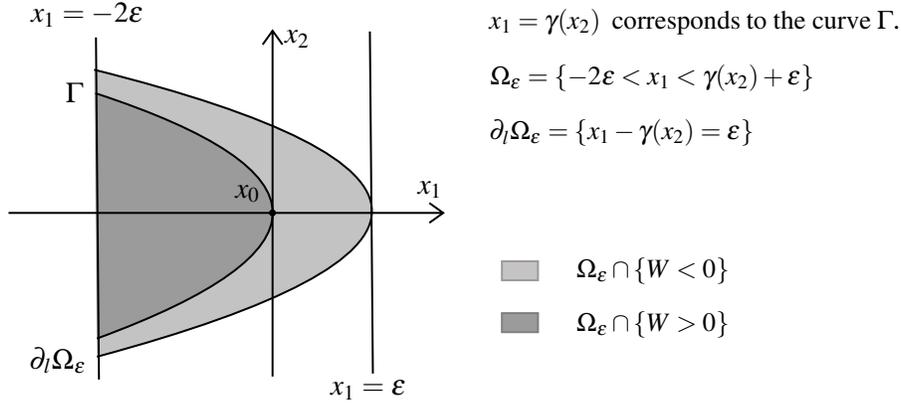
\draw (200,273.16) node [anchor=north west][inner sep=0.75pt]  [font=\small] [align=left] {{\fontfamily{ptm}\selectfont \textbf{Fig. 1.} \ Moving plans around the nodal set.}};
\draw (132.63,106) node   [align=left] {\begin{minipage}[lt]{10.27pt}\setlength\topsep{0pt}
$\Gamma $
\end{minipage}};
\draw (106,233) node [anchor=north west][inner sep=0.75pt]   [align=left] {$\partial _{l} \Omega_{\varepsilon }$};
\draw (106.93,59.14) node [anchor=north west][inner sep=0.75pt]   [align=left] {$x_{1} =-2\varepsilon $};
\draw (258.07,249.69) node [anchor=north west][inner sep=0.75pt]   [align=left] {$x_{1} =\varepsilon $};
\draw (302,149) node [anchor=north west][inner sep=0.75pt]   [align=left] {$x_{1}$};
\draw (235,73) node [anchor=north west][inner sep=0.75pt]   [align=left] {$x_{2}$};
\draw (210,151) node [anchor=north west][inner sep=0.75pt]   [align=left] {$x_{0}$};
\draw (338.44,90.4) node [anchor=north west][inner sep=0.75pt]  [font=\small] [align=left] {$\Omega_{\varepsilon } =\{-2\varepsilon < x_{1} < \gamma (x_{2} )+\varepsilon \}$};
\draw (382,188) node [anchor=north west][inner sep=0.75pt]  [font=\small] [align=left] {$\Omega_{\varepsilon } \cap \{W< 0\}$};
\draw (382.00,215) node [anchor=north west][inner sep=0.75pt]  [font=\small] [align=left] {$\Omega_{\varepsilon } \cap \{W >0\}$};
\draw (224,159) node [anchor=north west][inner sep=0.75pt]  [font=\Huge] [align=left] {$\cdot $};
\draw (338,62) node [anchor=north west][inner sep=0.75pt]   [align=left] {{\small $x_{1} =\gamma (x_{2} )${\fontfamily{ptm}\selectfont \, corresponds to the curve }$\Gamma $}.};
\draw (338.44,131.4) node [anchor=north west][inner sep=0.75pt]  [font=\small] [align=left] {$ $};
\draw (338,118) node [anchor=north west][inner sep=0.75pt]  [font=\small] [align=left] {$\partial _{l} \Omega_{\varepsilon } =\{x_{1} -\gamma (x_{2} )=\varepsilon \}$};

\end{tikzpicture}

\end{center}

\bigskip

First, let us mention the key estimate when the singular source is positive. Recall the auxiliary function $v=u-w+C_0(\varepsilon+\gamma(x_2)-x_1)$ with certain harmonic function $w$ and its equation (see (3.16) and (3.17) in \cite{jmpa-4}):
\begin{equation}\label{equ-moving}
    \Delta v+f(x,v(x))-C_0\gamma''(x_2)=0,\quad \mbox{in}\ \ \Omega_{\varepsilon},
\end{equation}
where the function $f(x,v)$ now becomes
$$f(x,v)=|x|^{2\alpha}W(x)e^{v+w-C_0(\epsilon+\gamma(x_2)-x_1)}.$$
In order to make the method of moving planes feasible, we need to show
\begin{equation}\label{main-ine}
    f(x,v)\le f(x_{\lambda},v), \ \ \mbox{for every } x\in\Omega_{\varepsilon}\cap\{x_1\ge\lambda\}   \ \ \mbox{with} \; -\varepsilon_1<\lambda<\varepsilon,
\end{equation}
which is described as
\begin{equation}\label{pos-com}
|x|^{2\alpha}W(x)e^{v+w-C_0(\epsilon+\gamma(x_2)-x_1)}<|x_{\lambda}|^{2\alpha}W(x_{\lambda})e^{v+w_{\lambda}-C_0(\epsilon+\gamma(x_2)-x_1^{\lambda})}.
\end{equation}
Here $x_{\lambda}=(2\lambda-x_1,x_2)$ and $\varepsilon_1\in(0,\varepsilon)$ is certain small constant.

\smallskip

\noindent{\bf Proof of (\ref{pos-com}):} In the first situation, if $W(x^{\lambda})>0$ and $W(x)\le 0$, there is nothing to prove. The second situation is when both $W(x)$ and $W(x^{\lambda})$ are positive. In this case, we observe that $\lambda<0$. Consequently, we have
$|2\lambda-x_1|>|x_1|$, which implies $|x_{\lambda}|^{2\alpha}>|x|^{2\alpha}$. Because of this, all other terms can be estimated as in \cite{jmpa-4}. In other words, the new term is helpful.

The final case is when both $W(x)$ and $W(x^{\lambda})$ are negative, which means $x$ and $x_{\lambda}$ are both in the region $\{W<0\}$. In this case, if we ignore the function $w$ because the $C_0$ part can easily majorize it, we need to have 
$$|x|^{2\alpha}W(x)<|x_{\lambda}|^{2\alpha}W(x^{\lambda})e^{2C_0(\lambda-x_1)}.$$
If we write $W(x)$ as 
$$W(x)=-d(x)h(x)$$
where $d(x)$ is the distance to the curve $x_1=\gamma(x_2)$ from the region $W<0$, $h$ is a smooth function bounded above and below by two positive constants. Since the difference between $h(x)$ and $h(x_{\lambda})$ is also majorized by the $C_0$ part, we ignore it. The inequality we prove in this case is 
\begin{equation}\label{key-1a}
2\alpha \log |x|+\log d(x)\ge 2\alpha \log |x_{\lambda}|+\log d(x_{\lambda}). 
\end{equation}
This clearly reduces to 
\begin{equation}\label{key-1b}
\log \Big(1+\frac{d(x)-d(x_{\lambda})}{d(x_{\lambda})}\Big)\ge \alpha \log\Big(1+ \frac{4 |\lambda |\cdot |\lambda-x_1 |}{|x|^2}\Big).
\end{equation}
It is well known that $d(x)-d(x_{\lambda})\ge C|x_1-\lambda |$. 
Since the curve $\Gamma$ is $x_1=\gamma(x_2)$ and $\Gamma$ is convex, we see that $|x_1|\le a x_2^2$ for some $a>0$. Thus for $|x|$ small we have $|x|>>|x_1|$. To compare the terms, we also observe that $|x|$ and $|x_{\lambda}|$ are comparable. Therefore when we choose $\varepsilon>0$ small enough such that $|\lambda |$ small, we see that (\ref{key-1b}) holds. So (\ref{key-1a}) is justified. Consequently, we complete the proof of \eqref{pos-com} in all situations.

\qed

\medskip

\noindent{\textbf{Remark 2.1.} We would like to point out that the study of blow-up solutions for the mean field equation has been greatly advanced by the recent outstanding works of Wei-Zhang \cite{wei-zhang-adv,wei-zhang-plms,wei-zhang-jems}. We expect our method to be used for other equations as well (see \cite{ahmedou-wu-zhang}). }

\end{document}